\documentclass{gtart_a}
\pdfoutput=1

\usepackage{pinlabel,stmaryrd}

\title[Dynamics of the mapping class group action on character varieties]{Dynamics of the mapping class group action\\on the variety of
$\mathrm{PSL}_2 \mathbb{C}$ characters}

\author{Juan Souto}
\givenname{Juan}
\surname{Souto}
\address{Department of Mathematics\\
University of Chicago\\\newline
5734 University Avenue\\
Chicago IL 60637-1514\\
USA}
\email{juan@math.uchicago.edu}
\urladdr{}

\author{Peter Storm}
\givenname{Peter}
\surname{Storm}
\address{Department of Mathematics\\
Stanford University\\\newline
450 Serra Mall\\
Stanford CA 94305-2125\\
USA}
\email{pastorm@stanford.edu}
\urladdr{}

\volumenumber{10}
\issuenumber{}
\publicationyear{2006}
\papernumber{20}
\lognumber{0702}
\startpage{715}
\endpage{736}

\doi{}
\MR{}
\Zbl{}


\received{14 January 2006}
\revised{16 May 2006}
\accepted{30 April 2006}
\published{11 July 2006}
\publishedonline{11 July 2006}
\proposed{Benson Farb}
\seconded{Walter Neumann, Jean-Pierre Otal}
\corresponding{}
\editor{}
\version{}

\def\notin{\not\in}
\def\co{\mskip 0.5mu\colon\thinspace}

\makeatletter
\def\cnewtheorem#1[#2]#3{\newtheorem{#1}{#3}[section]
\expandafter\let\csname c@#1\endcsname\c@sat}

\newtheorem{sat}{Theorem}[section]
\cnewtheorem{lem}[sat]{Lemma}
\cnewtheorem{kor}[sat]{Corollary}
\cnewtheorem{prop}[sat]{Proposition}
\theoremstyle{definition}
\cnewtheorem{bei}[sat]{Example}
\cnewtheorem{defi}[sat]{Definition}
\makeatother  

\makeautorefname{sat}{Theorem}
\makeautorefname{kor}{Corollary}
\makeautorefname{bei}{Example}

\newtheorem*{namedtheorem}{\theoremname}
\newcommand{\theoremname}{testing}
\newenvironment{named}[1]{\renewcommand{\theoremname}{#1}\begin{namedtheorem}}{\end{namedtheorem}}

\theoremstyle{remark}

\newcommand{\BC}{\mathbb C}			\newcommand{\BH}{\mathbb
H}
\newcommand{\BR}{\mathbb R}			
\newcommand{\BN}{\mathbb N}			
			\newcommand{\BZ}{\mathbb
Z}

\newcommand{\CC}{\mathcal C}

\newcommand{\CO}{\mathcal O}		\newcommand{\CP}{\mathcal P}
\newcommand{\CQ}{\mathcal Q}

		\newcommand{\CX}{\mathcal X}

\renewcommand{\D}{\partial}

\DeclareMathOperator{\Out}{Out}		
\DeclareMathOperator{\Diff}{Diff}	
\DeclareMathOperator{\SL}{SL}		
\DeclareMathOperator{\PSL}{PSL}		
\DeclareMathOperator{\Hom}{Hom}		
\DeclareMathOperator{\vol}{vol}		

\DeclareMathOperator{\SU}{SU}
\DeclareMathOperator{\Map}{\mathrm{Mod}}

\begin{document}

\begin{asciiabstract}
We study the action of the mapping class group Mod(S) on the boundary dQ
of quasifuchsian space Q. Among other results, Mod(S) is shown to be
topologically transitive on the subset C in dQ of manifolds without a
conformally compact end.  We also prove that any open subset of the
character variety X(pi_1(S),SL(2,C)) intersecting dQ does not admit a
nonconstant Mod(S)-invariant meromorphic function. This is related to a
question of Goldman.
\end{asciiabstract}

\begin{htmlabstract}
We study the action of the mapping class group Mod(S) on the
boundary &part;Q of quasifuchsian space Q. Among
other results, Mod(S) is shown to be topologically transitive
on the subset C &sub; &part;Q of manifolds without
a conformally compact end.  We also prove that any open subset of the
character variety X(&pi;<sub>1</sub>(S),SL<sub>2</sub><b>C</b>)
intersecting &part;Q does not admit a nonconstant
Mod(S)&ndash;invariant meromorphic function. This is related to
a question of Goldman.
\end{htmlabstract}

\begin{webabstract}
We study the action of the mapping class group $\mathrm{Mod}(S)$ on the
boundary $\partial\mathcal{Q}$ of quasifuchsian space $\mathcal{Q}$. Among
other results, $\mathrm{Mod}(S)$ is shown to be topologically transitive
on the subset $\mathcal{C}\subset\partial\mathcal{Q}$ of manifolds without
a conformally compact end.  We also prove that any open subset of the
character variety $\mathcal{X}(\pi_1 (S), \mathrm{SL}_2 \mathbb{C})$
intersecting $\partial\mathcal{Q}$ does not admit a nonconstant
$\mathrm{Mod}(S)$--invariant meromorphic function. This is related to
a question of Goldman.
\end{webabstract}

\begin{abstract}
We study the action of the mapping class group $\Map(S)$ on the boundary
$\D\CQ$ of quasifuchsian space $\CQ$. Among other results, $\Map(S)$ is
shown to be topologically transitive on the subset $\CC\subset\D\CQ$ of
manifolds without a conformally compact end.  We also prove that any open
subset of the character variety $\CX (\pi_1 (S), \SL_2 \BC)$ intersecting
$\D \CQ$ does not admit a nonconstant $\Map(S)$--invariant meromorphic
function. This is related to a question of Goldman.
\end{abstract}

\maketitle

\section{Introduction}

Let $S$ be a closed oriented surface of genus $g\ge 2$ and let $\Gamma$
be its fundamental group. The mapping class group
$$\Map(S)=\Diff(S)/\Diff_0(S)=\Out(\Gamma)$$
of $S$ acts on the character variety
$$\CX(\Gamma,\PSL_2\BC)=\Hom(\Gamma,\PSL_2\BC) \sslash \PSL_2\BC$$
by precomposition. Quasifuchsian space $\CQ$ is the open cell in the
character variety $\CX(\Gamma,\PSL_2\BC)$ formed by the conjugacy classes
of faithful representations with convex cocompact image.  It is invariant
under the mapping class group, and the action of $\Map(S)$ on $\CQ$
is properly discontinuous. On the other hand, our first result shows
that the action of $\Map(S)$ on $\D\CQ$ has very complicated dynamics.

\begin{sat}\label{topo}
Let $\CC \subset \wwbar{\CQ}$ denote the set of representations whose
quotient manifold has no conformally compact end and let
$\overline{\CC}$ denote the closure of $\CC$.  Then:
\begin{enumerate}
\item $\overline{\CC}$ is a $\Map(S)$--invariant
nowhere dense topologically perfect set.
\item The action of $\Map(S)$ on $\overline{\CC}$ is topologically
transitive.
\item The points $\rho\in\D\CQ$ satisfying $\overline{\CC}
\subset\overline{\Map(S) \cdot \rho}$ form a dense $G_\delta$--set.
\end{enumerate}
\end{sat}

In particular, \fullref{topo} implies that any continuous
$\Map(S)$--invariant function on $\D\CQ$ is constant.  Recall that the
action of a group on a locally compact Hausdorff separable topological
space is topologically transitive if the translates of any two open sets
intersect, or equivalently, if there is a dense orbit.

The group $\PSL_2\BC$ is the group of orientation preserving isometries
of hyperbolic 3-space $\BH^3$.	It is well known that $\rho$ is
faithful and that the action of $\rho(\Gamma)$ on $\BH^3$ is free and
properly discontinuous for each $\rho\in\wwbar\CQ$ (see Kapovich
\cite[Theorem~9.1.4]{Kap}).  In particular $M_\rho=\BH^3/\rho(\Gamma)$
is an orientable hyperbolic manifold homotopy equivalent to $S$. From
this point of view the set $\CC$ of \fullref{topo} is the set of all
$\rho\in\wwbar\CQ$ such that the boundary of the convex core of the
associated hyperbolic manifold $M_\rho$ does not have a compact component.

A special role is played by geometrically finite representations
$\rho \in \CC$ where the boundary of the convex core of $M_\rho$ is a
collection of thrice-punctured spheres.  We call these representations
full maximal cusps.  The proof of \fullref{topo} involves studying the
dynamics of the $\Map(S)$--action near these points.  For example, the
techniques used to prove \fullref{topo} are also used to prove that for
any open neighborhood $U \subset \D \CQ$ of a full maximal cusp, the orbit
$\Map(S) \cdot U$ is dense in $\D \CQ$ (see \fullref{cor of main prop 2}).

It is well known that two representations $\rho,\rho'\in\CQ$ are close if
and only if the associated hyperbolic manifolds $M_\rho$ and $M_{\rho'}$
are bi-Lipschitz with a small bi-Lipschitz constant. This is why geometric
invariants of the manifold $M_\rho$, for example the volume of the convex
core, the injectivity radius, the lowest eigenvalue of the Laplacian,
or the Hausdorff dimension of the limit set, are continuous functions
on quasifuchsian space.  However, in the larger set $\wwbar{\CQ}$
the picture is more complex.  Two representations $\rho, \rho' \in
\wwbar{\CQ}$ may be close without there being any bi-Lipschitz
homeomorphism from $M_\rho$ to $M_{\rho'}$.  With this motivation one may
ask which geometric invariants remain continuous on $\wwbar{\CQ}$.
In each case it was previously known via different methods that these
quantities are no longer continuous on $\wwbar\CQ$.	We derive from
\fullref{topo} a unified proof of this fact:

\begin{sat}\label{klein}
The volume of the convex core, the injectivity radius, the lowest
eigenvalue of the Laplacian and the Hausdorff dimension of the limit
set do not vary continuously on $\wwbar{\CQ}$.
\end{sat}

We also apply \fullref{topo} to study $\Map(S)$--invariant
meromorphic functions defined on subsets of the character variety
$\CX(\Gamma,\SL_2\BC)$. In \cite{Goldman03}, Goldman proved that
every $\Map(S)$--invariant meromorphic function defined on the whole
of $\CX(\Gamma, \SL_2\BC)$ must be constant. This result motivates the
question of which connected open subsets $U$ of $\CX(\Gamma,\SL_2\BC)$
admit nonconstant $\Map(S)$--invariant meromorphic functions. (A weaker
form of this question can be found in \cite[Section~1.4]{Goldman03}.)
Goldman \cite{Goldman97} deduced the nonexistence of invariant nonconstant
meromorphic functions from the ergodicity of the $\Map(S)$--action on the
(real) subvariety $\CX(\Gamma,\SU_2)$.  In particular his result applies
to every connected open set $U\subset\CX(\Gamma,\SL_2\BC)$ containing
unitary representations. We obtain the following analogue of Goldman's
result:

\begin{sat}\label{Goldman question}
Let $U \subset \CX (\Gamma, \SL_2 \BC)$ be a $\Map(S)$--invariant connected
open set. If $U$ contains both (faithful) convex cocompact representations
and indiscrete representations then any $\Map(S)$--invariant meromorphic
function on $U$ is a constant function.
\end{sat}

After some preliminaries in \fullref{prelims} we present in \fullref{main
construction} a concrete construction of hyperbolic 3-manifolds which
will be one of the key ingredients in the proof of \fullref{topo}.
Full maximal cusps are the other key ingredient. In particular, in
\fullref{sec:ubiquity} we deduce that full maximal cusps are dense in
$\CC$.	The proof uses techniques developed by McMullen \cite{McMullen},
Canary, Culler, Hersonsky, and Shalen \cite{CCHS,Canary-Hersonsky},
who studied the ubiquity of maximally cusped representations on the
boundaries of various deformation spaces.  \fullref{klein} is proved
in \fullref{sec:klein} and \fullref{Goldman question} is proved in
\fullref{sec:goldman}.

\rk{Acknowledgements}
The first author is supported by a postdoctoral fellowship of the Deutsche
Forschungs Gemeinschaft.  The second author is supported by a National
Science Foundation Postdoctoral Fellowship.
The first author would like to thank the members of the University of
Chicago Department of Mathematics for their hospitality.  The second
author would like to thank Richard Canary and Misha Kapovich for helpful
conversations on these topics.

\section{Preliminaries } \label{prelims}

We refer to Heusener and Porti \cite{Porti-Heusener} for basic
facts about the character variety and to Anderson \cite{Anderson} for
a survey about the deformation theory of discrete subgroups of
$\PSL_2\BC$. If $H$ is a finitely generated non-virtually abelian
torsion free group then $\Hom(H,\PSL_2\BC)$ is a complex algebraic
variety on which $\PSL_2\BC$ acts by conjugacy. The character variety
$$\CX(H,\PSL_2\BC)=\Hom(H,\PSL_2\BC) \sslash \PSL_2\BC$$
is the quotient of $\Hom(H,\PSL_2\BC)$ under this action
in the sense of invariant theory. We remind the reader that
$\CX(H,\PSL_2\BC)$ does not coincide with the set theoretic quotient
$\Hom(H,\PSL_2\BC)/\PSL_2\BC$. However, the set of conjugacy classes of
discrete faithful representations is contained in a smooth open manifold
in $\CX(H,\PSL_2\BC)$ \cite[Section 4]{Porti-Heusener}.  This paper is
concerned only with discrete faithful representations, so the machinery
of invariant theory will not be needed.

\begin{named}{Notation}
The Greek letters $\rho$ and $\sigma$ (possibly with decoration)
will be used to indicate \emph{conjugacy classes} of representations.
Thus $\rho$ and $\sigma$ will be elements of the appropriate character
variety. The notation $\rho(H)$ will indicate the image of $H$ under
any fixed homomorphism of the conjugacy class $\rho$.
\end{named}

We identify the group $\PSL_2\BC$ with the group of orientation preserving
isometries of hyperbolic 3-space $\BH^3$.  We will denote the convex core
of a hyperbolic manifold $M$ by $CC(M)$.  The convergence of sequences
in $\CX (H,\PSL_2\BC)$ is said to be \emph{algebraic convergence}.
Let $\rho_i\to\rho$ be an algebraically convergent sequence.  If the
representations $\rho_i$ are discrete and faithful for all $i$,
then it is well known that $\rho$ is discrete and faithful as well
(see Kapovich \cite[Theorem~9.1.4]{Kap}). Moreover, up to a choice
of a subsequence, the groups $\rho_i(H)$ converge in the Chabauty
topology to a discrete subgroup $H_G$ of $\PSL_2\BC$ which contains the
image of $\rho$. $H_G$ is the \emph{geometric limit} of the sequence
$\{\rho_i(H)\}$.

In this paper we will mainly consider discrete and faithful
representations of the fundamental group $\Gamma$ of a closed surface
$S$. A faithful and discrete representation $\rho\in\CX(\Gamma,\PSL_2\BC)$
induces a homotopy class of homotopy equivalences $S \longrightarrow
M_{\rho}$.  (Recall that $M_\rho$ denotes the hyperbolic $3$--manifold
$\BH^3 / \rho(\Gamma)$.)  A theorem of Bonahon \cite{Bonahon86}
ensures that $M_{\rho}$ is homeomorphic to a trivial interval bundle
over $S$.  Moreover, the homotopy equivalence $S\longrightarrow M_{\rho}$
determines a unique isotopy class of orientation preserving homeomorphisms
$S\times(-1,1) \longrightarrow M_{\rho}$.  In other words, $M_\rho$
has a positive end (the top end) and a negative end (the bottom end).
Observe that orientation reversing elements in $\Map(S)$ extend in a
canonical way to orientation preserving homeomorphisms of $S\times(-1,1)$
which interchange the top and bottom ends.

A component of $\partial CC(M_\rho)$ is said to face the top
(resp. bottom) end of $M_\rho$ if it is isotopic in $M_\rho - CC(M_\rho)$
out the top (resp. bottom) end of $M_\rho$.  With this terminology,
the top (resp. bottom) end of $M_{\rho}$ is \emph{conformally compact}
if there is a single compact component of $\partial CC(M_\rho)$ facing
the top (resp. bottom) end of $M_{\rho}$, or equivalently if there are
compact embedded convex surfaces exiting the top (resp. bottom) end of
$M_{\rho}$.  (This terminology comes from the fact that a conformally
compact end limits onto a compact boundary component of the conformal
manifold $ \left( \BH^3 \cup \Omega_\rho \right) / \rho(\Gamma)$,
where $\Omega_\rho \subset \mathbb{S}^2_\infty$ denotes the domain of
discontinuity of $\rho(\Gamma)$.)

Quasifuchsian space $\CQ \subset \CX(\Gamma, \PSL_2 \BC)$ is the open
set of conjugacy classes of faithful convex cocompact representations.
The closure $\wwbar\CQ$ of quasifuchsian space $\CQ$ consists of
faithful representations with discrete image and hence it is contained
in an open submanifold of $\CX(\Gamma, \PSL_2 \BC)$. Let $\partial \CQ$
denote the boundary of quasifuchsian space $\wwbar\CQ - \CQ \subset
\CX(\Gamma, \PSL_2 \BC)$.  Sullivan \cite{Sullivan} proved that $\D\CQ$
is also the set
$$\bigl( \overline{\CX(\Gamma, \PSL_2 \BC) - \wwbar\CQ } \bigr) -
\left( \CX(\Gamma, \PSL_2 \BC) - \wwbar\CQ \right).$$
In other words, $\partial \CQ$ is the frontier of $\wwbar{\CQ}$.
Note that $\rho\in\wwbar\CQ$ is quasifuchsian if and only if it has
two conformally compact ends.

The following $\Map(S)$--invariant subset of $\D\CQ$ will play a central
role:
$$\CC  \ = \{\rho\in\D\CQ \ \vert \ M_\rho\ \hbox{has no conformally
compact end}\}.$$
For example, if $CC(M_\rho) = M_\rho$ then $\rho \in \CC$.  In particular,
a cyclic cover of a fibered hyperbolic $3$--manifold lies in $\CC$.
On the other hand, the closure of any Bers slice is disjoint from $\CC$.
Finally, any representation in $\wwbar{\CQ}$ without parabolic elements
has a neighborhood disjoint from $\CC$ (see the proof of \fullref{topo}).
This gives a decomposition of $\wwbar{\CQ}$ into three pieces:
the manifolds with two conformally compact ends $\CQ$, the manifolds
with exactly one conformally compact end $\D \CQ - \CC$, and finally
$\CC$, which contains several types of manifolds.  Roughly speaking,
the $\Map(S)$--action becomes increasingly chaotic on these pieces.
Surprisingly, $\CC$ is not closed (see Section 7).


\section{The main construction} \label{main construction}
In this section we present a construction which is the core of the
proof of \fullref{topo}. The main building pieces in our construction
are so called maximal cusps.  Fix a compact hyperbolic surface $S$
with fundamental group $\Gamma$.

We will say that a discrete finitely generated subgroup of $\PSL_2\BC$
is a \emph{full maximal cusp} if it is geometrically finite and every
component of the boundary of the convex core of the associated hyperbolic
manifold is a thrice punctured sphere.	Observe that $\rho\in\CC$ if
$\rho(\Gamma)$ is a full maximal cusp.	We will say a representation
$\rho \in \D\CQ$ is a \emph{one sided maximal cusp} if it is geometrically
finite, has one conformally compact end, and each component of $\partial
CC(M_\rho)$ facing the end of $M_\rho$ which is not conformally compact
is a thrice punctured sphere. (A one sided maximal cusp is often simply
called a maximal cusp (see McMullen \cite{McMullen}).	Our modified
terminology has been chosen for clarity.)  The set of full maximal cusps
in $\D \CQ$ is countable, and intuitively forms a set of ``rational
points" on the boundary.  This intuition can be made precise in the
punctured torus case.

Maximal cusps are very convenient when making concrete constructions
because any pair of totally geodesic thrice punctured spheres in any pair
of hyperbolic 3-manifolds are isometric. In particular, if $M$ and $M'$
are hyperbolic manifolds whose convex cores boundaries $\D CC(M)$ and $\D
CC(M')$ contain thrice punctured spheres $X$ and $X'$, and $\phi\co
X\longrightarrow X'$ is a homeomorphism, then $\phi$ is isotopic to an
isometry which we denote again by $\phi$.  Hence there is a hyperbolic
manifold $N$ which is covered by $M$ and $M'$, whose convex core is
isometric to $CC(M)\cup_\phi CC(M')$.

The second main ingredient in our constructions is the following
lemma, which is a consequence of Thurston's Dehn-filling Theorem
(see Thurston \cite[Chapter~4]{Thurston}, Bonahon and Otal
\cite{Bonahon-Otal}, and Comar \cite{Comar}).

\begin{lem}\label{Dehn}
Let $H < \PSL_2\BC$ be a geometrically finite group such that $\BH^3/H$
is homeomorphic to $S\times(0,1){-}\CP$ where $\CP$ is an unlinked
collection of disjoint simple closed curves. Then the group $H$ is
the geometric limit of geometrically finite groups $H_n$ isomorphic to
$\pi_1(S)$. Moreover, if $H$ is a full (resp. one sided) maximal cusp,
then the $H_n$ can be chosen to be full (resp. one sided) maximal
cusps.
\end{lem}

Recall that a collection of disjoint simple closed curves in
$S\times[0,1]$ is unlinked if every curve is contained in an embedded
boundary parallel surface which is disjoint from all the other curves. The
groups $H_n$ in the statement of \fullref{Dehn} are obtained by performing
a hyperbolic $(1,n)$--Dehn surgery on a neighborhood of each of the curves
in $\CP$.

\begin{lem}\label{lem:ag}
Let $H, H_n < \PSL_2 \BC$ be as in the statement of \fullref{Dehn}, and
let $\rho_n\co \Gamma\longrightarrow H_n$ be isomorphisms.  If $\rho\co
\Gamma \longrightarrow \PSL_2 \BC$ is the representation induced by an
embedded level surface
$$S \times \{ t \} \subset \left( S \times (0,1) - \CP \right) \cong
\BH^3/H$$
then there are automorphisms $\alpha_n$ of $\Gamma$ such that $\rho_n
\circ \alpha_n$ converges algebraically to $\rho$.
\end{lem}
\begin{proof}
Since $H_n \to H$ geometrically, for any a compact submanifold $K \subset
\BH^3/H$ there is a sequence of smooth embeddings
$\phi_n\co K \longrightarrow \BH^3/H_n$ which converge in
the $\mathcal{C}^\infty$--topology to isometric embeddings
(see McMullen \cite[Section~2.2]{McM2}).  Since $\BH^3/H_n$ is obtained
by $(1,n)$--Dehn surgery on $\BH^3/H$, it is clear that by choosing a
sufficiently large compact submanifold $K$ the restricted homomorphism
$$(\phi_{n*}) |_{\rho(\Gamma)}\co \rho(\Gamma) \longrightarrow H_n <
\PSL_2 \BC$$
will be well defined and injective, and thus an isomorphism.  Since the
maps $\phi_n$ are converging to isometries, the sequence of homomorphisms
$(\phi_{n*})|_{\rho(\Gamma)}$ is converging to the identity map.
Therefore
$$\phi_{n*} \circ \rho \co \Gamma \longrightarrow H_n$$
is a sequence of discrete faithful representations converging to $\rho$.
Define $\alpha_n$ to be $\rho_n^{-1} \circ \phi_{n*} \circ \rho$.
This proves the lemma.
\end{proof}

The following is the main result of the present section.  For clarity
it has been split into three similar pieces.


\begin{prop} \label{main prop a}
Let $\rho,\bar\rho\in\D\CQ$ be one sided maximal cusps.  There exists
a sequence of representations $\{\rho_i\}$ in quasifuchsian space and
a sequence of mapping classes $\{\alpha_i \}$ such that
$$\rho_i \to \rho \quad \text{and} \quad ( \alpha_i \cdot \rho_i )
\to \bar\rho.$$
\end{prop}

\begin{prop} \label{main prop b}
Let $\rho \in \D\CQ$ be a one sided maximal cusp.  Let $\rho^{\CC}$ be a
full maximal cusp.  There exists a sequence of one sided maximal cusps $\{
\rho_i \}$ and a sequence of mapping classes $\{ \alpha_i \}$ such that
$$\rho_i \to \rho \quad \text{and} \quad ( \alpha_i \cdot \rho_i )
\to \rho^{\CC}.$$
\end{prop}

\begin{prop} \label{main prop c}
Let $\rho^\CC, \bar\rho^\CC \in \D\CQ$ be full maximal cusps.  There
exists a sequence of full maximal cusps $\{ \rho_i \}$ and a sequence
of mapping classes $\{ \alpha_i \}$ such that
$$\rho_i \to \rho^\CC \quad \text{and} \quad ( \alpha_i \cdot \rho_i )
\to \bar\rho^{\CC}.$$
\end{prop}

The proofs of Propositions \ref{main prop a}, \ref{main prop b}, and
\ref{main prop c} are very similar, so we will prove carefully only the
first one.

\begin{proof}
Up to composition with an orientation reversing element in $\Map(S)$,
we may assume that the top end of $M_\rho$ and the bottom end of
$M_{\bar\rho}$ are conformally compact.

By \fullref{lem:ag} it suffices to construct a hyperbolic manifold
$N$ homeomorphic to the complement in $S\times(0,1)$ of an unlinked
collection of simple curves, which is (locally isometrically) covered
by both $M_\rho$ and $M_{\bar\rho}$.

\begin{figure}[p]
\begin{center}
\labellist\small
\pinlabel {$CC(M_\rho)$} [l] at 293 567
\pinlabel {$M$} at 318 398
\pinlabel {$\CP_+$} [l] at 368 502
\pinlabel {$\phi_+(\CP_+)$} [l] at 368 442
\pinlabel {$\phi_-(\CP_-)$} [l] at 368 356
\pinlabel {$\CP_-$} [l] at 368 294
\pinlabel {$CC(M_{\bar\rho})$} [l] at 293 235
\pinlabel {glue by $\phi_+$} [l] at 318 470
\pinlabel {glue by $\phi_-$} [l] at 318 320
\endlabellist
\includegraphics[scale=0.9]{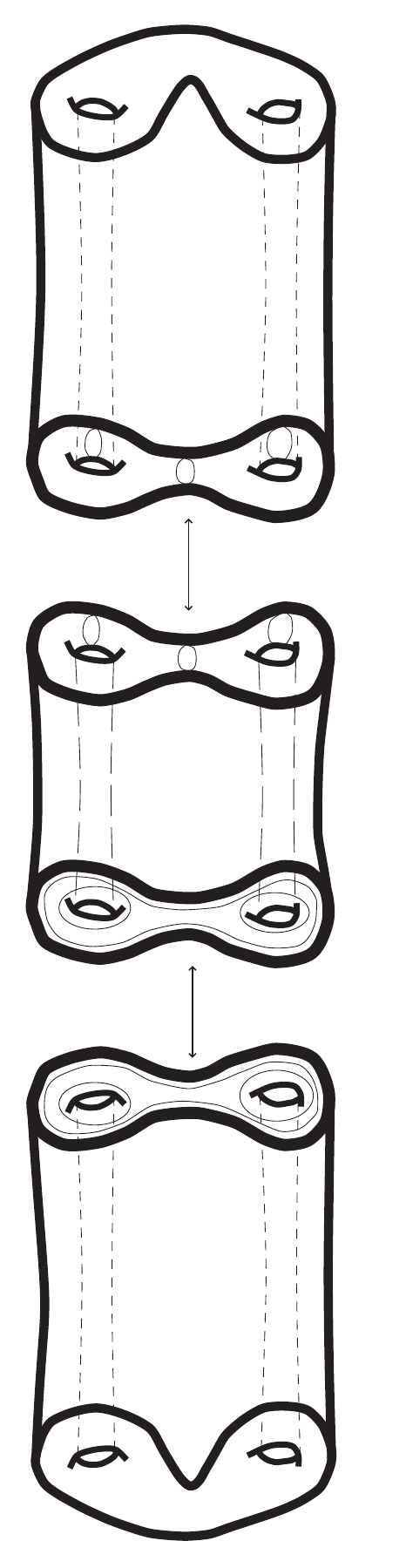}
\caption{Constructing the manifold $N$} \label{figure1}
\end{center}
\end{figure}
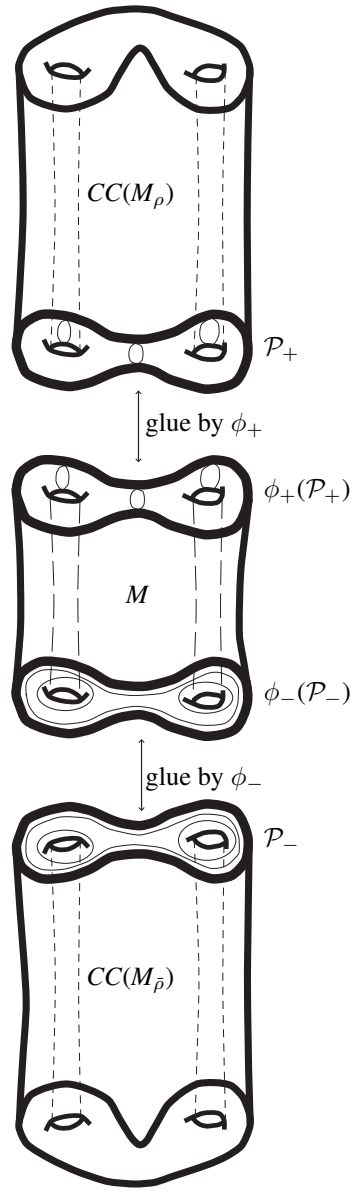

The convex cores $CC(M_{\rho})$ and $CC(M_{\bar\rho})$ are by definition
homeomorphic to
$$\left( S\times[-1,1] \right) \  - \  (\CP_+ \times\{-1\})\ \hbox{and}\
\left(S\times[-1,1] \right) \ - \ (\CP_- \times\{1\}),$$
where $\CP_+$ and $\CP_-$ are pants decompositions of $S$. Choose
sufficiently complicated homeomorphisms $\phi_+\co S\longrightarrow S$ and
$\phi_- \co S\longrightarrow S$ such that the two collections
$\phi_+(\CP_+)$ and $\phi_-(\CP_-)$ bind the surface $S$. Thurston's
hyperbolization theorem (see Otal \cite{Otal-Haken}) implies that there is a
geometrically finite hyperbolic 3-manifold $S\longrightarrow M$ whose
convex core is homeomorphic to
$$(S\times[-1,1])\ {-} \
    \left( \phi_-(\CP_-)\times\{-1\}\cup\phi_+(\CP_+)\times\{1\}
    \right).$$
The middle manifold $M$ is a full maximal cusp.  $M_{\rho}$ and
$M_{\bar\rho}$ are one sided maximal cusps.  Therefore the convex cores
of $M_{\rho}$, $M_{\bar\rho}$ and $M$ can be glued together to obtain the
convex core $CC(N)$ of a hyperbolic 3-manifold $N$ covered by $M_{\rho}$
and $M_{\bar\rho}$ and homeomorphic to
$$(S\times[-3,3]) \ {-}\
   \left( \phi_- (\CP_-)\times\{-1\}\cup\phi_+(\CP_+)\times\{1\}
   \right).$$
(See \fullref{figure1} for a picture of this construction.  In the figure
the homeomorphisms $\phi_{\pm}$ are the identity map, because the pants
decompositions $\CP_{\pm}$ already bind $S$.)

We are now in the situation of \fullref{Dehn}.	The holonomy
representation of $N$ corresponds to $H < \PSL_2 \BC$, and
$$\phi_- (\CP_-)\times\{-1\} \ \cup \ \phi_+(\CP_+)\times\{1\}$$
corresponds to the unlinked collection of curves $\CP$.  Two applications
of \fullref{lem:ag} conclude the proof of \fullref{main prop a}.
First apply \fullref{lem:ag} to the level surface $S \times \{ 2 \}$.
Then apply it to the level surface $S \times \{ -2 \}$.
\end{proof}

For the proof of \fullref{main prop b}, replace $\bar\rho$ with $\rho^\CC$
and follow the above argument.	This will require adding a third pants
decomposition to the above notation, corresponding to the bottom end
of $M_{\rho^\CC}$.  Similarly, for the proof of \fullref{main prop c},
replace $\rho$ with $\rho^\CC$, $\bar\rho$ with $\bar\rho^\CC$, and
follow the above argument.  In this case it will be necessary to add a
fourth pants decomposition to the above notation.  Otherwise the proof
goes through unchanged in both cases.


\section{Ubiquity of maximal cusps}\label{sec:ubiquity}

In this section we study the ubiquity of maximal cusps in $\D\CQ$.
In particular, full maximal cusps are shown to be dense in the subset
$\CC \subset \D\CQ$ given by manifolds without a conformally compact
end (see \fullref{cusps}).  This result will be used in the proof
of \fullref{topo}.  We will use techniques developed by McMullen
\cite{McMullen}, Canary, Culler, Hersonsky and Shalen \cite{CCHS}, and
Canary and Hersonsky \cite{Canary-Hersonsky}.  We begin by combining results
from \cite{McMullen} and \cite{CCHS} with work of Brock, Bromberg,
Evans and Souto \cite{BBES} to prove the following statement.

\begin{sat} \label{McMullen thm}
One sided maximal cusps form a dense subset of $\D\CQ$.
\end{sat}

\begin{proof}
A standard Baire category argument shows that the set of representations
$\rho \in \D\CQ$ without parabolic elements forms a dense subset of
$\D\CQ$ (see McMullen \cite[Corollary~1.5]{McMullen} or Canary, Culler,
Hersonsky and Shalen \cite[Lemma~15.2]{CCHS}).	So it suffices to
approximate only representations without parabolic elements by one sided
maximal cusps.	Pick such a representation $\rho \in \D\CQ$.  There are
only two possibilities: either the limit set of $\rho$ is the entire
sphere at infinity, or $M_\rho$ has exactly one conformally compact end
and one geometrically infinite end.

If the limit set of such a $\rho$ is the entire sphere at infinity, then
Canary, Culler, Hersonsky and Shalen have proven \cite[Theorem~6.1]{CCHS}
that $\rho$ is an algebraic limit of full maximal cusps.  \fullref{main
prop b} shows that a full maximal cusp is always an algebraic limit of
one sided maximal cusps in $\D\CQ - \CC$.  Therefore it only remains to
consider the case where $M_\rho$ has one conformally compact end.

Let us assume without loss of generality that the top end of $M_\rho$ is
conformally compact.  Applying \cite{BBES}, we know the representation
$\rho \in \D\CQ$ is a strong limit of quasifuchsian representations
$\rho_i$.  By Kerckhoff and Thurston \cite{Kerckhoff-Thurston}, the
conformal structure at infinity of the top end of $M_\rho$ must be the
limit in Teichm{\"u}ller space of the conformal structures at infinity
of the top ends of the $M_{\rho_i}$.  Therefore altering the manifolds
$M_{\rho_i}$ by a $K_i$--bi-Lipschitz deformation, where $K_i \rightarrow
1$, produces a sequence of quasifuchsian manifolds in a Bers slice
converging strongly to $M_\rho$.  We may now apply \cite{McMullen} to
conclude that $\rho$ is an algebraic limit of one sided maximal cusps.
\end{proof}

From \fullref{main prop b} and \fullref{McMullen thm} we derive the
following ostensibly stronger result:

\begin{prop} \label{main prop 2}
In every open subset of $\D\CQ$ we find a $\bar\rho$ with the following
property:  For every full maximal cusp $\rho^\CC$ there is a sequence
$\{\alpha_i\}$ in $\Map(S)$ with
$$\rho^\CC =\lim_i\alpha_i \cdot \bar\rho$$
\end{prop}

\begin{proof}
Recall that the set of all full maximal cusps is countable. In particular
there is a sequence of full maximal cusps $\{\rho_j\} \subset \CC$ such
that every full maximal cusp appears infinitely often in the sequence.

Moreover, we can fix open neighborhoods $V_j$ of $\rho_j$ in $\D\CQ$ such
that when a subsequence $\rho_{j_i}$ converges then for every choice of
$\rho_{j_i}'\in V_{j_i}$ the sequence $\rho_{j_i}'$ also converges to
the same limit.

Given an open neighborhood $U_1$ in $\D\CQ$ we obtain from \fullref{main
prop b} and \fullref{McMullen thm} a one sided maximal cusp $\sigma_1\in
U_1 \subset \D\CQ$ and $\alpha_1\in\Map(S)$ with $\alpha_1 \cdot
\sigma_1 \in V_1$.  Let $U_2$ be an open and relatively compact
neighborhood of $\sigma_1$ in $U\cap\alpha_1^{-1}(V_1)$. Applying again
the same argument we find a one sided maximal cusp $\sigma_2 \in U_2$
and $\alpha_2 \in \Map(S)$ with $\alpha_2 \cdot \sigma_2 \in V_2$.
Let $U_3$ be an open and relatively compact neighborhood of $\sigma_2$ in
$U_2\cap\alpha_2^{-1}(V_2)$. Inductively we find a sequence of one sided
maximal cusps $\{\sigma_i\}$, a sequence $\{\alpha_i\}$ in $\Map(S)$,
and a sequence of open sets $\{U_i\}$ such that: $U_{i+1}$ is relatively
compact in $U_i$, $\sigma_i$ is in $U_i$, and $\alpha_i(U_i)\subset
V_i.$  By our choice of the neighborhoods $V_i$ we can conclude that
every $\bar\rho$ in $\cap_i U_i$ has the desired property.
\end{proof}

By reversing the logic of \fullref{main prop 2} we obtain:

\begin{kor} \label{cor of main prop 2}
If $W \subset \D \CQ$ is an open set containing a full maximal cusp then\break
$\Map(S)\cdot W$ is a dense open subset of $\D \CQ$.  Moreover, the set
$$\{ \rho \in \D \CQ \, | \, \overline{\Map(S)\cdot \rho} \text{  contains
every full maximal cusp} \}$$
is a dense $G_\delta$--set.
\end{kor}
\begin{proof}
The first sentence of the corollary follows immediately from \fullref{main
prop 2}.  To prove the second, let us reuse the notation $\{ \rho_j \}$
and $\{ V_j \}$ from the proof of \fullref{main prop 2}.  Then
$$\bigcap_j \Map(S) \cdot V_j = \{ \rho \in \D \CQ \, | \,
\overline{\Map(S) \cdot \rho} \text{  contains every full maximal cusp}
\}.$$
From this the corollary follows.
\end{proof}

The following theorem paraphrases the results of
Canary and Hersonsky \cite[Section~10]{Canary-Hersonsky}.  It will be
used in the proof of \fullref{cusps}.

\begin{sat}{\rm{\cite[Theorem~10.1]{Canary-Hersonsky}}}\label{CH theorem}\qua
Let $M$ be a hyperbolic $3$--manifold with a holonomy representation
$\rho$.  Assume there is a sequence of geometrically finite
representations $\rho_i$ converging algebraically to $\rho$ such that
the homomorphisms
$\rho (\pi_1(M)) \longrightarrow \rho_i (\pi_1(M))$
are induced by homeomorphisms
$$(M, \text{cusps of } M) \longrightarrow (M_{\rho_i}, \text{cusps of }
M_{\rho_i}).$$
Let $\mathcal{E}$ denote the geometrically infinite ends of
$M-\text{\{cusps\}}$.  Then there exists a sequence of geometrically
finite representations $\widehat{\rho}_j$ converging algebraically to
$\rho$ satisfying:
\begin{enumerate}
\item There exists a homeomorphism $\phi_j\co M \longrightarrow
M_{\widehat{\rho}_j}$ taking the cusps of $M$ into the cusps of
$M_{\widehat{\rho}_j}$.
\item If $\Sigma$ is a component of $\partial CC(M_{\widehat{\rho}_j})$
isotopic in
$$M_{\widehat{\rho}_j} - \text{int}(CC(M_{\widehat{\rho}_j}))$$
outside every bounded subset of $\phi_j (\mathcal{E})$, then $\Sigma$
is a thrice punctured sphere.
\end{enumerate}
\end{sat}

Using some standard terminology we have not defined here (see
Anderson \cite{Anderson}), this theorem can be restated informally as follows.
If $M$ lies on the boundary of the deformation space of a pared manifold,
then it can be approximated by geometrically finite manifolds on the
boundary of the same deformation space, where each geometrically infinite
end of $M - \text{\{cusps\}}$ has been replaced by a maximally cusped
geometrically finite end.  The statement of \fullref{CH theorem} differs
slightly from the statement of \cite[Theorem~10.1]{Canary-Hersonsky}.  However,
the proof of \cite[Theorem~10.1]{Canary-Hersonsky} proves \fullref{CH theorem}.

Recall that $\CC \subset \D\CQ$ is the set of manifolds without a
conformally compact end.  We are now ready to prove the main result of
this section, which will be used in the proof of \fullref{topo}.

\begin{prop} \label{cusps}
Full maximal cusps are dense in $\CC$.
\end{prop}

\begin{proof}
Pick a representation $\rho \in \CC$.  The goal is to approximate $\rho$
by a sequence of full maximal cusps.  We divide the proof into four cases,
depending on what type of representation $\rho$ is.

\textbf{Case 1}\qua Assume $\rho$ is geometrically finite.

We will produce a hyperbolic manifold covered by $M_\rho$ which is the
geometric limit of a sequence of full maximal cusps in $\CC$.  Then we
will apply \fullref{lem:ag} to complete Case 1.

The convex core $CC(M_\rho)$ of $M_\rho$ is homeomorphic to
$$N_0 := (S\times[-1,1]) {-} \left( \CP_-\times\{-1\}\cup\CP_+\times\{1\}
\right),$$
where $\CP_\pm$ are nonempty collections of disjoint, nonparallel simple
closed curves. Define
$$N := (S\times[-1,2]) \ {-} \	(\CP_-\times\{-1\}\cup\CP_+\times\{1\}).$$
We consider the embedding $N_0 \longrightarrow N$ and the induced morphism
$$\Pi\co\CX(\pi_1(N),\PSL_2\BC)\longrightarrow\CX(\Gamma,\PSL_2\BC).$$
By a result of Brock, Bromberg, Evans and Souto \cite[Lemma~5.1]{BBES}
there is an open set in $\Pi^{-1}(\rho)$ consisting of discrete, faithful,
geometrically finite, and minimally parabolic representations.  (Here,
minimally parabolic means that every parabolic conjugacy class of $\pi_1
(N)$ corresponds to a component of one of the pants decompositions
$\CP_-$ or $\CP_+$.)  Moreover, the marked conformal structure at
infinity corresponding to the bottom end of any manifold in the fiber
$\Pi^{-1}(\rho)$ is equivalent to the marked conformal structure at
infinity corresponding to the bottom end of $M_\rho$.

On the other hand, by deforming the added generator in any $\mathbb{Z}
\times \mathbb{Z}$ subgroup of $\pi_1(N)$, we see that every connected
component of $\Pi^{-1}(\rho)$ contains indiscrete representations.
So there must be a nonempty separating set of discrete faithful
representations in $\Pi^{-1}(\rho)$ which are either geometrically
infinite or not minimally parabolic.  A generic $\sigma \in\Pi^{-1}(\rho)$
is minimally parabolic.  (This is again a Baire category argument;
compare with McMullen \cite[Corollary~1.5]{McMullen} or Canary, Culler,
Hersonsky and Shalen \cite[Lemma~15.2]{CCHS}.)	We deduce
that there is a discrete, faithful, and minimally parabolic $\sigma \in
\Pi^{-1}(\rho)$ such that the top end of $M_{\sigma}$ is geometrically
infinite.

One should imagine that $M_\sigma$ has been obtained by adding a
geometrically infinite cap onto the top end of $M_\rho$.  With its
geometrically infinite end, $M_\sigma$ is more easily approximated by
one sided maximal cusps.  Indeed, we may now apply \fullref{CH theorem}
to conclude that $\sigma$ is the algebraic limit of a sequence $\sigma_i$
in $\CX(\pi_1(N),\PSL_2\BC)$ such $\sigma_i$ is discrete, faithful,
geometrically finite, and the top end of $M_{\sigma_i}$ is maximally
cusped.


Since $M_{\sigma}$ is minimally parabolic, we may apply a theorem of Evans
\cite{Evans} to conclude that the manifolds $M_{\sigma_i}$ also converge
geometrically to $M_{\sigma}$. \fullref{Dehn} implies now that for all
$i$ the manifold $M_{\sigma_i}$ is the geometric limit of a sequence
$M_i^j$ of geometrically finite hyperbolic manifolds homeomorphic to
$S\times(0,1)$ such that the top end of each $M_i^j$ is maximally cusped.
Taking a diagonal sequence we deduce that $M_{\sigma}$ is the geometric
limit of a sequence of geometrically finite manifolds whose top ends
are maximally cusped.  Since $M_{\sigma}$ is covered by $M_\rho$, we
may apply \fullref{lem:ag} to conclude that

\emph{$M_\rho$ is an algebraic limit of geometrically finite manifolds
in $\CC$ whose top ends are maximally cusped.}

Therefore, we may without loss of generality assume the the top end
of $M_\rho$ is maximally cusped.  Now apply the above argument again
to $M_\rho$, swapping the roles of the bottom end and the top end.
This yields a sequence of full maximal cusps converging algebraically
to $M_\rho$.  The completes the proof of Case 1.

\textbf{Case 2}\qua $M_\rho$ has at least one rank one cusp in each end.

By Case 1 it suffices to approximate $M_\rho$ by geometrically finite
manifolds in $\CC$. Let $P_+$ (resp. $P_-$) be a nonempty collection
of disjoint essential annuli on $S$ such that there is a relative
homeomophism
$$\left( S \times (0,1) \ , \ (P_+ \times [.8,1) \right) \bigcup \left(
P_- \times (0,.2]) \right) \longrightarrow(M_\rho, \text{cusps of }
M_\rho).$$
By a result of Brock, Bromberg, Evans and Souto
\cite[Corollary~3.2]{BBES}, $M_\rho$ is the algebraic and geometric limit
of a sequence of geometrically finite manifolds $M_{\rho_i}$ such that for
each $i$ there is a relative \emph{homotopy equivalence} $$\left( S \times
(0,1) \ , \ (P_+ \times [.8,1) \right) \bigcup \left(P_- \times (0,.2])
\right) \longrightarrow(M_{\rho_i}, \text{cusps of } M_{\rho_i}).$$
The concurrence of algebraic and geometric convergence implies that,
after passing to a subsequence, these relative homotopy equivalences
are rel homotopic to relative homeomorphisms.  This implies that the
manifolds $M_{\rho_i}$ are in $\CC$.

\textbf{Case 3}\qua $M_\rho$ has all of its rank one cusps in exactly
one end.

Assume without loss of generality that $M_\rho$ has rank one cusps in its
bottom end.  Then the top end of $M_\rho$ is necessarily geometrically
infinite.  The proof begins as in Case 3.  Let $P_-$ be a nonempty
collection of disjoint essential annuli on $S$ such that there is a
relative homeomophism
$$\left( S \times (0,1), P_- \times (0,.2] \right) \longrightarrow \left(
M_\rho, \text{cusps of } M_\rho \right).$$
As above, applying \cite[Corollary~3.2]{BBES} yields that $M_\rho$ is the
algebraic and geometric limit of a sequence of geometrically finite
manifolds $M_{\rho_i}$ such that for each $i$ there is a relative
homotopy equivalence
$$\left(S \times (0,1), P_- \times (0,.2] \right) \longrightarrow \left(
M_{\rho_i}, \text{cusps of } M_{\rho_i} \right).$$
Again using the fact that $M_{\rho_i}$ converges both algebraically
and geometrically to $M_\rho$, after possibly passing to a subsequence
we may assume these relative homotopy equivalences are homotopic to
relative homeomorphisms.

Having verified its hypotheses, we may now apply \fullref{CH theorem} to
find an $M_\rho$ convergent sequence of geometrically finite manifolds
in $\CC$ whose top ends are maximally cusped.  This reduces the proof
to Case 1.


\textbf{Case 4}\qua $M_\rho$ has no cusps.

In this case $M_\rho$ must have an empty domain of discontinuity.  We may
then directly apply \cite[Theorem~10.1]{Canary-Hersonsky} to conclude
that $M_\rho$ is the algebraic limit of a sequence of full maximal cusps.
\end{proof}

\section{The main theorem}\label{sec:klein}

In this section we prove the main result of this paper, \fullref{topo}.
As an application we also prove \fullref{klein}.

Recall that $\CC$ denotes the set of representations whose quotient
manifolds have no conformally compact end.  A goal of \fullref{topo}
is to prove that the closure $\overline{\CC}$ of $\CC$ is a small set,
but up to this point we have not found even a single discrete faithful
representation outside of $\overline{\CC}$.  This we now do, using a
theorem of Evans \cite{Evans}.

\begin{lem}\label{Evans lem}
If $\sigma \in \D \CQ$ is a representation with no parabolic elements
and exactly one conformally compact end (i.e. a singly degenerate
representation without parabolics), then $\sigma$ is not contained in
$\overline{\CC}$.
\end{lem}

\begin{proof}
Pick a representation $\sigma \in \D\CQ$ with no parabolic elements
and exactly one conformally compact end.  Suppose that $\sigma
\in \overline{\CC}$.  Then by \fullref{cusps} there is a sequence
of full maximal cusps $\tau_i \rightarrow \sigma$.  Since $\sigma$
has no parabolic elements, it follows from a theorem of Evans \cite{Evans}
that the manifolds $M_{\tau_i}$ converge geometrically to $M_\sigma$.
Since the top end of $M_\sigma$ is conformally compact, there is a
strictly convex surface $\Sigma$ embedded in the top end of $M_\sigma$.
Use the almost isometric embeddings provided by geometric convergence to
push the surface $\Sigma$ into $M_{\tau_i}$.  For sufficiently large $i$
this yields a strictly convex embedded surface in $M_{\tau_i}$, showing
that the manifolds $M_{\tau_i}$ eventually have a conformally compact
end.  Since they are all full maximal cusps, this is a contradiction.
Therefore $\sigma$ is not in $\overline{\CC}$.
\end{proof}


\begin{named}{\fullref{topo}}
Let $\CC \subset \CQ$ denote the set of representations whose quotient
manifold has no conformally compact end.  Then:
\begin{enumerate}
\item The closure $\overline{\CC}$ of $\CC$ is a $\Map(S)$--invariant
nowhere dense topologically perfect set.
\item The action of $\Map(S)$ on $\overline{\CC}$ is topologically
transitive.
\item The points $\rho\in\D\CQ$ satisfying $\overline{\CC}
\subset\overline{\Map(S) \cdot \rho}$ form a dense $G_\delta$--set.
\end{enumerate}
\end{named}

\begin{proof}
The set $\CC$ is $\Map(S)$--invariant by definition, implying the same for
its closure $\overline{\CC}$.  The set $\CC$ is topologically perfect
since full maximal cusps are dense in $\CC$ and \fullref{main prop c}
implies full maximal cusps are not isolated.

To finish claim (1), it remains to prove that $\overline{\CC}$ is nowhere
dense.	Following \fullref{Evans lem}, let $\sigma \notin \overline{\CC}$
be a representation without parabolics whose quotient manifold has
exactly one conformally compact end.  Let $U$ be a neighborhood of
$\sigma$ in the complement of $\overline{\CC}$.  By \fullref{main prop 2}
there is a representation in $U$ whose orbit limits onto all of $\CC$.
As $\Map(S) \cdot U \cap \overline{\CC} = \emptyset$, this proves that
$\overline{\CC}$ is nowhere dense.

We now prove claim (2).  Let $U^1$ and $U^2$ be open sets in $\CC$. By
\fullref{cusps} there are full maximal cusps $\rho^1 \in U^1$ and $\rho^2
\in U^2$. By \fullref{main prop c} there exist sequences $\{ \rho_i \}$
in $\CC$ and $ \{ \alpha_i \}$ in $\Map(S)$ such that
$$\rho^1 =\lim_i \rho_i \ \ \hbox{and}\ \  \rho^2 = \lim_i\alpha_i
\cdot \rho_i$$
This shows that for sufficiently large $i$, $(\alpha_i \cdot U^2) \cap
U^1$ is not empty. Hence the action of $\Map(S)$ on $\CC$ is topologically
transitive.  This implies topological transitivity on $\overline{\CC}$.

Finally, claim (3) follows immediately from \fullref{cusps} and
\fullref{cor of main prop 2}.
\end{proof}

As an application we prove the following theorem.

\begin{named}{\fullref{klein}}
The volume of the convex core, the injectivity radius, the lowest
eigenvalue of the Laplacian and the Hausdorff dimension of the limit
set do not vary continuously on $\wwbar{\CQ}$.
\end{named}

Each part of \fullref{klein} was known previously.  We present merely
a unified proof.

\begin{proof}
All these invariants are $\Map(S)$--invariant functions
$$\D\CQ\longrightarrow\BR\cup\{\infty\}.$$
In particular we derive from \fullref{topo} (3) and (4) that they either
are constant or fail to be continuous. In particular, it suffices to
find $\rho,\rho'\in\D\CQ$ for which they do not take the same value. Let
$\rho$ be a full maximal cusp and $\rho'$ the cyclic cover of the mapping
torus of a pseudo-Anosov $\phi\in\Map(S)$; $\rho'$ is in the boundary
of quasifuchsian space by the work of Thurston \cite{ThurstonII}.
We know that
$$\vol(CC(M_\rho))<\infty\ \hbox{and}\ \vol(CC(M_{\rho'}))=\infty$$
because $\rho$ is geometrically finite and $\rho'$ is not. The manifold
$M_{\rho'}$ covers a compact manifold and hence has positive injectivity
radius while $M_\rho$ has cusps and hence its injectivity radius vanishes.
Geometric finiteness of $\rho$ implies that $\Lambda_\rho$ has Hausdorff
dimension less than 2 (see Tukia \cite{Tukia} and Sullivan
\cite{Sullivan2}) while $\Lambda_{\rho'}=\BC P^1$ has dimension 2.

Sullivan's theorem tells us that the lowest eigenvalue of the Laplacian
on $M_\rho$ is
$$\text{dim} (\Lambda_\rho) \, ( 2 - \text{dim} (\Lambda_\rho)) > 0.$$
Finally, the lowest eigenvalue of the Laplacian on $M_{\rho'}$ equals
the infimum
$$ \inf_{f \in C^\infty_0 (M_{\rho'})}	\frac{ \int | \nabla f |^2 }{\int
| f |^2 }.$$
(This is true for any Riemannian manifold.)  Since $M_{\rho'}$ is the
cyclic cover of the mapping torus of $\phi$, it is easy to show without
any machinery that this infimum is zero.
\end{proof}

From this we see that it is in fact impossible to define a nontrivial
purely geometric invariant which varies continuously on $\wwbar{\CQ}$.

\section{$\text{Map}(S)$--invariant meromorphic
functions}\label{sec:goldman}

In this section we prove \fullref{Goldman question}.  As a warm up to
the $\SL_2\BC$ Case, we first consider meromorphic functions on open
subsets of $\CX(\Gamma,\PSL_2\BC)$.

\begin{sat}\label{warm-up thm}
Let $U \subset\CX(\Gamma,\PSL_2\BC)$ be a $\Map(S)$--invariant connected
open set. If $U$ contains both (faithful) convex cocompact representations
and indiscrete representations then any $\Map(S)$--invariant meromorphic
function on $U$ is a constant function.
\end{sat}

Note that indiscrete representations are dense in the complement of
quasifuchsian space (see Sullivan \cite{Sullivan}).

\begin{proof}
Without loss of generality we may assume that $U$ is a manifold. The set
$U$ intersects the interior and the exterior of $\wwbar\CQ$ and hence
$\D\CQ\cap U\neq\emptyset$. We claim that every continuous
$\Map(S)$--invariant function $f\co U\longrightarrow\BC$ is constant on
$\D\CQ\cap U$. Given one sided maximal cusps $\rho,\bar\rho\in\D\CQ \cap
U$ we obtain from \fullref{main prop a} a sequence $\{\rho_i\}\subset\CQ$
and a sequence of mapping classes $\{\alpha_i\}\subset\Map(S)$ with
$$\lim_i\rho_i=\rho\ \hbox{and}\ \lim_i\alpha_i \cdot \rho_i=\bar\rho.$$
We have $\rho_i\in U$ for sufficiently large $i$. The continuity and
the $\Map(S)$--invariance of $f$ imply that $f(\rho)=f(\bar\rho)$. Since
maximal cusps are dense in $\D\CQ\cap U$, the claim follows.

Assume now that $f\co U\longrightarrow\BC$ is meromorphic and
$\Map(S)$--invariant. The divisor $D$ of poles of $f$ is either empty or
has complex codimension $1$. The open set $U - D$ is open, connected,
$\Map(S)$--invariant and also intersects $\D\CQ$. In particular we may
assume without loss of generality that $f$ is holomorphic. We proved
above that $f$ is constant on the separating set $\D\CQ\cap U$. By
holomorphicity, this implies that $f$ is everywhere constant.
\end{proof}

We now discuss the relationship between the character varieties
$\CX(\Gamma,\SL_2\BC)$ and $\CX(\Gamma,\PSL_2\BC)$
and prove \fullref{Goldman question}.  The homomorphism
$\SL_2\BC\longrightarrow\PSL_2\BC$ induces a map
$$p\co \CX(\Gamma,\SL_2\BC)\longrightarrow\CX(\Gamma,\PSL_2\BC)$$
which maps onto the connected component of $\CX(\Gamma,\PSL_2\BC)$
containing $\wwbar{\CQ}$ (see Heusener and Porti
\cite[Remark~4.3]{Porti-Heusener} and Culler \cite{Culler}).
Recall that the closure of quasifuchsian space
$\wwbar{\CQ}\subset\CX(\Gamma,\PSL_2\BC)$ is contained in a smooth
open $\Map(S)$--invariant manifold $\CO\subset\CX(\Gamma,\PSL_2\BC)$
(see \cite[Section~4]{Porti-Heusener} and Goldman \cite{Goldman03}). After possibly shrinking
$\CO$ slightly we may assume that
$$p|_{p^{-1}(\CO)}\co p^{-1}(\CO)\longrightarrow\CO$$
is a Galois cover with covering transformation group $H_1(S,\BZ/2\BZ)$
(see \cite{Porti-Heusener}).
See Goldman \cite{Goldman03} for properties of $\CX(\Gamma,\SL_2\BC)$.

\begin{named}{\fullref{Goldman question}}
Let $U \subset\CX(\Gamma,\SL_2\BC)$ be a $\Map(S)$--invariant connected
open set. If $U$ contains both (faithful) convex cocompact representations
and indiscrete representations then any $\Map(S)$--invariant meromorphic
function on $U$ is a constant function.
\end{named}

Again note that indiscrete representations are dense in the complement
of the set of convex cocompact representations.

\begin{proof}
As above we may assume without loss of generality that $U$ is a manifold
contained in $\CO$. Moreover, it suffices to show that a holomorphic
$\Map(S)$--invariant function $f$ on $U$ is constant.

The conditions on $U$ imply that there is some $\sigma\in\D\CQ\cap p(U)$
since the image under $p$ of a faithful convex cocompact representation
lies in $\CQ$ and the image of an indiscrete representation is again
indiscrete. Choose a neighborhood $V'$ of $\sigma$ and $V\subset U$ open
such that $p\vert_{V}\co V\longrightarrow V'$ is a homeomorphism. Since
the actions of $\Map(S)$ on $\CX(\Gamma,\SL_2\BC)$ and
$\CX(\Gamma,\PSL_2\BC)$ are compatible we obtain that the restriction
of $p$ to the open $\Map(S)$--invariant set
$$W= \bigcup_{\alpha\in\Map(S)} \alpha \cdot V$$
is a covering. It suffices to show that the restriction of $f$ to $W$
is constant. Given $n\in\BN$ we consider the holomorphic function on
$p(W)$ given by
$$S_n(z)=\sum_{w \in p^{-1}(z)} f^n(w).$$
We claim that $S_n$ is constant for all $n$. The function $S_n$ is
holomorphic and $\Map(S)$--invariant, but we cannot directly apply
\fullref{warm-up thm} since $p(W)$ may not be connected. However,
every connected component of $p(W)$ intersects $\D\CQ$ and therefore
the same proof as in \fullref{warm-up thm} applies to show that $S_n$
is locally constant.

The following lemma shows that $S_n$ being locally constant for all $n$
implies that $f$ is itself constant on $W$.

\begin{lem} \label{juan's a clever guy}
Let $\Omega \subset\BC^r$ be a connected open set and
$f_1,\dots,f_k\co \Omega \longrightarrow \BC$ holomorphic such that for
all $n$ the
function $S_n(z)=\sum f_i(z)^n$ is constant. Then all the functions $f_i$
are constant.
\end{lem}
\begin{proof}
Before going further observe that it suffices to consider the case where
$\Omega$ is a unit disk in $\BC$, because a nonconstant holomorphic
function will have a nonconstant restriction to some holomorphic disk.
Seeking a contradiction assume that the lemma is false.  One can
assume without loss of generality that none of the functions $f_i$
are constant and that no two of the functions $f_i$ and $f_j$ are
proportional. We may further assume, up to reducing $\Omega$, that for
all $i$ and $z$ we have $f_i(z)\neq 0$ and $f'_i(z)\neq 0$. Moreover,
the assumption that no two of the functions are proportional
implies that, up to relabeling, there is some $z_0$ satisfying
$|f_1(z_0)|>|f_2(z_0)|\ge\dots\ge|f_k(z_0)|$. Multiplying by a suitable
scalar we may assume that $f_1(z_0)=1$. Computing the derivative of $S_n$
at the point $0$ we obtain the following identity for all $n$:
$$0=S'_n(z_0)=n(f'_1(z_0)+f_2(z_0)^{n-1}f'_2(z_0)+\dots+f_k(z_0)^{n-1}f'_k(z_0))$$
Dividing by $n$ and taking a limit $n\to\infty$ we derive that $f'_1
(z_0)=0$. This is a contradiction.
\end{proof}

As mentioned above this concludes the proof of \fullref{Goldman question}.
\end{proof}

\section{$\CC \subset \D\CQ$ is not closed} \label{sad story}

This section will show that the set
$$\CC := \{ \rho \in \D\CQ \, | \, M_\rho  \text{  has no conformally
compact end} \}$$
(see \fullref{prelims}) is not a closed subset of $\D\CQ$.  This fact
surprised the authors.	This section is logically independent of Sections
\ref{main construction}--\ref{sec:goldman}.

The proof is a slight elaboration of a construction due to McMullen
\cite[Lemma~A.4]{McM3} (which was in turn based on results of Kerckhoff
and Thurston \cite{Kerckhoff-Thurston} and Anderson and Canary \cite{AC}).
As Lemma A.4 is very clearly written, we will not attempt to reproduce
its construction here.	We will assume that the reader has read Lemma A.4
and the example which follows it.  (These total only two and half pages.)
Our notation will be chosen to conform to McMullen's.

Choose two pants decompositions $P_{\pm}$ which bind the surface $S$.
Choose an embedded essential closed curve $C \subset S$ such that:
$C$ and $P_+$ bind $S$, and $C$ and $P_-$ bind $S$.  By Thurston's
hyperbolization theorem (see Otal \cite{Otal-Haken}) there is an
infinite volume hyperbolic $3$--manifold $N$ whose convex core is a finite
volume manifold with totally geodesic boundary homeomorphic to
$$(S \times [0,1] ) - \left( P_+ \times \{1\} \, \cup \, P_- \times \{0 \}
\, \cup \, C \times \{1/2 \} \right).$$
With this manifold $N$, perform McMullen's construction in
\cite[Lemma~A.4]{McM3}.  Following his notation, let $N_n$ be a sequence of
$(1,n)$--Dehn surgeries on $C \times \{ 1/2 \} \subset N$ together with
maps $F_n \co N \longrightarrow N_n$ converging in the
compact--$\mathcal{C}^\infty$ topology to an isometric embedding.  Let
$f\co S \longrightarrow N$ be an immersed essential surface which wraps
around $C$.

Mark the manifolds $N_n$ by the composition $F_n \circ f$.  Equipped
with these markings the sequence $N_n$ is contained in $\mathcal{C}
\subset \D\CQ$ and converges algebraically to the covering space of
$N$ given by $f_*(\pi(S))$.  Recall that the curve $C$ binds $S$ with
either pants decomposition $P_+$ or $P_-$.  From this it follows that
the only parabolics of $f_* (\pi_1 (S))$ correspond to the curve $C$.
Therefore the algebraic limit must have one conformally compact end,
and does not lie in $\CC \subset \D\CQ$.  This shows that $\CC$ is not
closed and concludes the construction.

It would be interesting to find a geometric characterization of the
manifolds in the closure $\overline{\CC}$.  This appears to be difficult.

\bibliographystyle{gtart}
\bibliography{link}

\end{document}